\documentclass{paper}

\usepackage[T1]{fontenc}
\usepackage[english]{babel}
\usepackage{amsmath}
\usepackage{amssymb}
\usepackage{amsthm}
\usepackage{natbib}
\usepackage{color}

\usepackage{amsthm}
\newtheorem{theo}{Theorem}
\newtheorem{lemm}[theo]{Lemma}

\usepackage{epsfig}

\def\O{{\mathcal{ O}}}

\newcommand{\F}{\textnormal{F}}

\def\R{\mathbb{R}}

\newcommand{\E}{\textnormal{E}}

\def\N{\mathbb{N}}

\renewcommand{\F}{{\cal F}}

\renewcommand{\F}{\mathcal{F}}

\newcommand{\eps}{\varepsilon}

\renewcommand{\P}{\mathcal{P}}

\newcommand{\Xn}{\mathbf{X^n}}

\DeclareMathAlphabet{\mathpzc}{OT1}{pzc}{m}{it}
\newcommand{\fha}{\hat{f}}

\newcommand{\I}{\mathbb{I}}

\renewcommand{\O}{\mathcal{O}}
\makeatother

\begin{document}

\title{Concentration rate and consistency of the posterior under monotonicity constraints}

\author{{Jean-Bernard} {Salomond}}

\institution{CREST - Universit\'e Paris Dauphine\\ 3 avenue Pierre Larousse \\ 92245 Malakoff, France\\}

\maketitle

\begin{abstract}
{In this paper, we consider the well known problem of estimating
a density function under qualitative assumptions. More precisely, we estimate monotone non increasing densities in a Bayesian
setting and derive concentration rate for the posterior distribution
for a Dirichlet process and finite mixture prior. We prove that the posterior distribution based on both priors
concentrates at the rate $(n/\log(n))^{-1/3}$, which is the minimax
rate of estimation up to a $\log(n)$ factor. We also study the behaviour of the posterior for the point-wise loss at any fixed point of the support the density and for the sup norm. We prove that the posterior is consistent for both losses.}
\end{abstract}

\paragraph*{keyword}

{Density estimation}
{Bayesian inference}
{Concentration Rate}

\section{Introduction}
\label{sec:intro}
The nonparametric problem of estimating monotone curves,
and monotone densities in particular, has been well studied in the
literature both from a theoretical and applied perspectives. Shape constrained estimation is
fairly popular in the nonparametric literature and widely used in practice
\citep[see][for instance]{robertson:wright:dykstra:1988}. Monotone densities appear in a wide variety of applications such as survival
analysis, where it is natural to assume that the uncensored survival time has a monotone non increasing density. In these problems, estimating the
survival function is equivalent to estimate the survival time density say $f$ and the pointwise estimate $f(0)$. It is thus interesting to have a better understanding of the behaviour of the estimation procedures in this case. An interesting property of monotone non increasing densities on $\R^+$ is
that they have a mixture representation pointed out by \citet{will:56} 
\begin{equation}\label{eq:Will}
f(x) = \int_{0}^\infty {\I_{[0,\theta]}(x) \over \theta } dP(\theta),
\end{equation}
where $P$ is a probability distribution on $\R^+$ called the mixing distribution. In order to emphasize the dependence in $P$, we will denote $f_P$
the functions admitting representation \eqref{eq:Will}. This representation allows for inference based on the likelihood. \citet{Gre1956} derived 
the nonparametric maximum likelihood estimator of a monotone density and \citet{PrakasaRao70}
studied the behavior of the Grenander estimator at a fixed point. \citet{MR822052} and more recently, \citet{balabdaoui:wellner} studied
very precisely the asymptotic properties of the non parametric maximum likelyhood estimator. It
is proved to be consistent and to converge at the minimax rate $n^{-1/3}$
when the support of the distribution is compact. In their paper \cite{Durot2012} get some refined asymptotic results for the supremum norm.   \\

The mixture representation of monotone densities lead naturally to a mixture type prior on the set of monotone non
increasing densities with support on $[0,L]$ or $\R^+$. 
For example \citet{MR736538} and \citet{lo84} introduced the Dirichlet Process prior (DP) and \citet{BrLo1989} considered the special case of unimodal densities with a prior based on a Dirichlet Process mixture. The problem of deriving concentration rates for mixtures
models have receive a huge interest in the past decade.
\citet{WG08} studied properties of general mixture models 
\citet{vdv2001} studied the well known problem of Gaussian mixtures, \citet{Rousseau10} derive
concentration rates for mixtures of betas, \citet{kruijer:rousseau:vdv:09} proved good adaptive properties of mixtures of Gaussian. Extensions to the multivariate case have recently been introduced (e.g. \citet{SGmulti12}). \\

Under monotonicity constrained, we derive an upper bound for the posterior concentration rate with respect to some metric or semi metric $d(\cdot,\cdot)$, that is a positive sequence $(\epsilon_n)_n$ that goes to $0$ when $n$ goes to infinity such that 
$$
\E_0^n\left ( \Pi(d(f,f_0) > \epsilon_n | \Xn) \right ) \to 0,
$$
where the expectation is taken under the true distribution $P_0$ of the data $\Xn$ and where $f_0$ is the density of $P_0$ with respect to the Lebesgue measure. 
 Following \citet{KHAZAEI:2010:INRIA-00494692:1} we study two families
of nonparametric priors on the class of monotone non increasing densities. Interestingly in our setting, the so called Kullback-Leibler property, that is the fact that the prior puts
enough mass on Kulback-Leibler neighbourhood of the true density, is not satisfied. Thus the
approach based on the seminal paper of \citet{ggv00} cannot be applied. We therefore use a modified version of their results and obtain for the two
families of prior a concentration rate of order $(n/\log(n))^{-1/3}$ which is the minimax estimation rate up to a $\log(n)$ factor under the $L_1$ or Hellinger distance.  We extend these
results to densities with support on $\R^+$ and prove that under some conditions
on the tail of the distribution, the posterior still concentrates at an almost optimal rate. To the author's knowledge, no concentration rates have been
derived for monotone densities on $\R^+$.

Interestingly, the non parametric maximum likelyhood estimator of $f_P(x)$
 is not consistent for $x = 0$ (see \citet{MR1392133} and \citet{balabdaoui:wellner}
 for instance). However, we prove that the posterior distribution of $f$
 is still consistent at this point under a specific family of non parametric mixture prior. In fact we prove the pointwise consistency of the posterior for all $x$ in $[0,L]$ with $L \leq \infty$. We then
 derive a consistent Bayesian estimator of the density at any fixed point of the support. This is
 particularly interesting as the point-wise loss is usually difficult to study in a Bayesian framework as the Bayesian approaches are well suited to
 losses related to the Kullback-Leiber divergence. We also study the behaviour of the posterior distribution for the sup norm when the density has a compact support. This problem has been
 addressed recently in the frequentist literature by \citet{Durot2012}. They derive refined asymptotic results on the sup norm of the difference
 between a Grenander-type estimator and the true density on sub intervals of the form $[\epsilon,L-\epsilon]$ where $\epsilon>0$ avoiding the problems at the boundaries. Here, we prove that the posterior distribution is consistent in
 sup norm on the whole support of $f_0$ when it has compact support. We also derive concentration rate for the posterior of the
 density taken at a fixed point and for the sup norm on subsets of $[0,L]$ for $L<\infty$. We also derive an upper bound for the concentration rate of $f(x)$ for $x \in (0,L)$ but only get suboptimal rates using a testing approach as in \citet{GineNickl2010}. It is to be noted that for this problem the modulus of continuity for the pointwise and Hellinger losses defined for $f_0 \in \F$ and $x \in (0,L)$ by 
 $$
 m(\epsilon) := \sup \{ |f(x) - f_0(x)|: f \in \F,~ h(f,f_0)\leq \epsilon\}  
 $$
 is of the order $\epsilon^{2/3}$ \citep[see][]{donoho1991geometrizing}. Given the discussion in \citet{Hebert2013concent}, it is to be expected that the usual approach of \citet{ggv00} based on tests will lead to suboptimal concentration rates.  
We now introduce some notations which will be needed throughout the paper. 
\paragraph{Notations}
For $0<L\leq \infty$ define the set $\F_L$ by
$$
\F_L = \left \{f\text{ s.t. } 0 \leq f < \infty ,\; f\searrow \, \;  \int_0^L f = 1\right\},
$$
We also define $\mathfrak{S}_k$ the $k$-simplex that is the set $\{(s_1, \dots , s_k )\in [0,1]^k , \sum_{i=1}^k s_i = 1\}$. Let $KL(p_1,p_2)$ be
the Kullback Leibler deviation between the densities $p_1$ and $p_2$ with respect to some measure $\lambda$ 
$$
KL(p_1,p_2) = \int \log\left( {p_1 \over p_2 } \right) p_1 d\lambda.  
$$
We also define the Hellinger distance $h(p_1,p_2)$ between $p_1$ and $p_2$ as 
$$
h^2(p_1,p_2) =\frac{1}{2} \int (\sqrt{p_1} - \sqrt{p_2} )^2 d\lambda .
$$
We will say that $\Xi^n = o_{p_0} (1)$ if $\Xi^n \to 0$ under $P_0$. Finally we will denote $f'$ the derivative of $f$. 

\paragraph{Construction of a prior distribution on $\F_L$}
Using the mixture representation of monotone non
increasing densities \eqref{eq:Will} we construct nonparametric priors on the
set $\F_L$ by considering a prior on the mixing distribution $P$. Let $\mathcal{P}$ be the set of probability measures on $[0,L]$. Thus we fall in
the well known set up of nonparametric mixture priors models. We consider two types of
prior on the set $\mathcal{P}$. 
\begin{description}
\item[Type 1 : Dirichlet Process prior] $P \sim DP(A,\alpha)$ where $A$ is a
positive constant and $\alpha$ a probability density on $[0,L]$. \\
\item[Type 2 : Finite mixture] $P = \sum_{j=1}^K p_j \delta_{x_j}$ with $K$ a
non zero integer and $\delta_x$ the dirac function on $x$. We choose a
prior distribution $Q$ on $K$ and given $K$, define distributions
$\pi_{x,K}$ on $(x_1, \dots, x_K) \in [0,L]^K$ and $\pi_{p,K}$ on $(p_1, \dots, p_K) \in \mathfrak{S}_K$.
\end{description}
For $\Xn = (X_1, \dots, X_n)$, a sample of $n$ independent and identically distributed
random variables with common probability distribution function $f$ in $\F_L$
with respect to the Lebesgue measure, we denote $\Pi(\cdot|\Xn)$ the posterior
probability measure associated with the prior $\Pi$. \\ 

The paper is organised as follow: the main results are given in Section
\ref{sec:main}, where conditions on the priors are discussed. The proofs are
presented in Section \ref{sec:proof}. 

\section{Main results} \label{sec:main}

Concentration rates of the posterior distributions have been well studied in the
literature and some general results link the rate to the prior (see
\cite{ggv00}). However, in our setting, the Kullback Leibler property is not satisfied in its usual form and thus the standard Theorems do not hold. In fact an interesting feature of mixture distributions whose kernels have varying support is that the prior mass of the sets $\{f, KL(f_0,f) = + \infty \}$ is $1$ for most $f_0 \in \F_L$ given that $f$ and $f_0$ will have different support. One could prevent this by imposing that the support of the mixing distribution is wider than the support of $f_0$, however this could lead to a deterioration of the concentration rate.  Here, we use a modified version of the results of \citet{ggv00} considering truncated versions of the density $f$. This idea has been considered in \citet{KHAZAEI:2010:INRIA-00494692:1} in a similar setting. 
We impose some conditions on the prior under which the posterior distribution concentrates at the minimax rate up to a $\log(n)$ term. 
\paragraph{Conditions on the prior} 
\begin{description}
\item[C1 condition on $\alpha$] Let $\alpha$ be a probability density on $\R^+$ such that for all $\theta \in (0,L)$, $\alpha(\theta)>0$. Consider the following conditions on $\alpha$ 
\begin{subequations}
\begin{itemize}
\item for $0< t_1 \leq t_2$ and $\theta$ small enough 
\begin{equation}
\theta^{t_1} \lesssim \alpha(\theta) \lesssim \theta^{t^2}
\label{eq:cond:alpha:a}
\end{equation}
\item for $1<a_1\leq a_2$ and $\theta$ small enough 
\begin{equation}
e^{-a_1/\theta} \lesssim \alpha(\theta) \lesssim e^{-a_2/\theta}
\label{eq:cond:alpha:b}
\end{equation}
\item for $1<b_1\leq b_2$ and $\theta$ small enough 
\begin{equation}
e^{-b_1/\theta} \lesssim \alpha(L-\theta) \lesssim e^{-b_2/\theta}
\label{eq:cond:alpha:c}
\end{equation}
\end{itemize}
\end{subequations}

\item[C2 condition for Type I prior] For $P \sim DP(\alpha,M)$ with $\alpha$ satisfying C1

\item[C3 condition for the Type II prior] The following conditions holds
\begin{itemize}
\item For some positive constants $C_1, C_2,a_1, \dots, a_k,c$ 
\begin{eqnarray}
 e^{-C_1 K \log(K)} \geq Q(K) \geq e^{-C_2 K \log(K)} \label{eq:cond_Q}
 \\
\pi_{p,k}(p_1,\dots, p_K) \geq K^{-K} c^K p_1^{a_1}\dots p_K^{a_K}
\label{eq:cond_p}
\end{eqnarray} 
\item $\pi_{x,K}$ is the distribution of $K$ independent and identically distributed random variables sampled from $\alpha$.

\end{itemize}
\item[C4 Condition for densities on $\R^+$] If $f_0 \in \F_\infty$ then for $\beta$ and $\tau$ some fixed positive constant we have for $x$ large enough 
\begin{equation}
f_0(x) \leq e^{-\beta x^\tau} .
\label{eq:cond:tails}
\end{equation}
\end{description}

\subsection{Posterior concentration rate for the $L_1$ and Hellinger metric}

The following Theorems gives the posterior concentration rate for the $L_1$ and Hellinger metric for monotone non increasing densities on $[0,L]$ with $L<\infty$ and $L = \infty$. For both Theorems the proofs are postponed to section \ref{sec:proof}.

\begin{theo}\label{th1}
Let $\Xn = (X_1,\dots,X_n)$ be an independent and identically distributed sample with a common probability
distribution function $f_0$ such that $f_0\in \F_L$ with $0<L<\infty$. Let $\Pi$ be either a Type I or Type II prior satisfying \textbf{C2} or \textbf{C3} respectively with $\alpha$ satisfying \eqref{eq:cond:alpha:a}. If $d(\cdot,\cdot)$ is either the $L^1$ or Hellinger distance, then  there exists a positive constant $C$ such that 
\begin{equation}
 \Pi\left(f,d(f,f_0)\geq C \left({n \over \log(n)} \right)^{-1/3}|\Xn\right)
\to 0, \qquad P_0 \, a.e. 
\label{eq:conc:helling:finite}
\end{equation}
when $n$ goes to infinity, where $C$ depends on $f_0$ only through $L$ and an upper bound on $f_0(0)$. Furthermore, if for $\delta >0$, $\sup_{[0,\delta]} |f_0'(x)| < \infty$ and $\alpha$ satisfies \eqref{eq:cond:alpha:b}, or $\sup_{[L,L-\delta]} |f_0'(x)| < \infty$ and $\alpha$ satisfies \eqref{eq:cond:alpha:c}, then \eqref{eq:conc:helling:finite} still holds.
\end{theo}
Conditions C1 and C2 are roughly the same as in \citet{KHAZAEI:2010:INRIA-00494692:1}. Theorem \ref{th1} is thus an extension of their results to concentration rates. We also extend their results to mixtures prior satisfying \eqref{eq:cond:alpha:b} or \eqref{eq:cond:alpha:c} under some additional conditions on $f_0$. This will prove useful for the estimation of $f_0$ and $f_L$. Under condition C3 on the tail of the true density, i.e.  we require exponential tails, we get the
posterior concentration rate for density with support on $\R^+$.
\begin{theo}
	Let $\Xn = (X_1,\dots,X_n)$ be an independent and identically distributed sample with a common probability distribution density $f_0$ such
	that $f_0 \in \F_\infty$ and $f_0$ satisfy \textbf{C3}. Let $\Pi$ be either a Type I or Type II prior satisfying \textbf{C2} or \textbf{C3} respectively with $\alpha$ satisfying \eqref{eq:cond:alpha:a}. Then for some positive constant $C$ we have for $d(\cdot,\cdot)$ either the $L_1$ or Hellinger metric
	\begin{equation}
		\Pi\left( d(f_P,f_0) \geq C \left( n/\log(n) \right)^{-1/3} \log(n)^{1/\tau} |\Xn \right) \to 0, \; P_0 \; a.e. 
		\label{eq:cvrate_infty}
	\end{equation}
	when $n$ goes to infinity. Similarly, if for $\delta >0$, $\sup_{[0,\delta]} |f_0'(x)| < \infty$ and $\alpha$ satisfies \eqref{eq:cond:alpha:b}, \eqref{eq:cvrate_infty} still holds. 
\label{th15}
\end{theo}
Note that considering monotone non increasing densities on $\R^+$ deteriorates the upper bound on the posterior concentration rate with a factor
$\log(n)^{1/\tau}$. It is not clear whether it could be sharpen or not. For instance, in the frequentist literature, \citet{RBRTM} observe a slower convergence rate when considering infinite support for densities without any other conditions. In a Bayesian setting, a similar log term appears in \citet{kruijer:rousseau:vdv:09} when considering densities with non compact support. However this deterioration of the concentration rate does not have a great influence on the asymptotic behaviour of the posterior. Note also
that the tail conditions are mild since $\tau$ can be taken as small as needed, and thus the considered densities can have almost polynomial tails. \\

The above results on  the posterior concentration rate in terms of the $L_1$ or Hellinger metric are new to our knowledge but not surprising. The specificity of these results lies in the fact that the usual approach based on the approach of \citet{ggv00} need to bound the prior mass of Kullback Leibler neighbourhoods of the true density which cannot be done here as explained in section \ref{sec:intro}. 

\subsection{Consistency and posterior concentration rate for the pointwise and supremum loss}
The following results consider the pointwise loss function for which only a few exist in the Bayesian nonparametric literature,  see for instance the paper of \citet{GineNickl2010}. The following Theorem proves consistency of the posterior distribution for all point in the interior of the support. 

\begin{theo}\label{th2}
Let $x$ be in $(0,L)$ with with $0<L\leq \infty$ but $x < \infty$. Let $f_0 \in \F_L$ such that $f_0'$ exists near $x$ and $f_0'(x) < 0$. Let $X_i$ , $i=1, \dots ,n $ and $\Pi$ be
either a Type I or Type II prior satisfying \textbf{C2} or \textbf{C3} respectively with $\alpha$ satisfying \textbf{C1} with either \eqref{eq:cond:alpha:a}, \eqref{eq:cond:alpha:b}  or \eqref{eq:cond:alpha:c}. Then, for all $x$ in $(0,L)$ with $x< \infty$, and $\epsilon >0$  
\begin{equation}
{\Pi}\big(|f_P(x) - f_0(x)|> \epsilon | \Xn \big) \to 0  .
\end{equation}
Consider the posterior median $\fha_n^\pi(x) = \inf \{ t, \Pi\big[f_P(x) \leq t | \Xn\big] > 1/2 \}$, it follows that  
\begin{equation}
P_0\big(|\fha^\pi_n(x) - f_0(x)| > \epsilon | \Xn \big) \to 0  .
\end{equation} 
\end{theo}
We thus have a pointwise consistency of the posterior distribution of $f_0(x)$ for every $x$ in the interior of the support of $f_0$. The maximum likelihood is not consistent at the boundaries of the support as pointed out in \citet{MR1392133} for
instance. In particular it is not consistent at $0$ and when $L< \infty$, it is not consistent at $L$. It is known that integrating the parameter as done in Bayesian approaches induces a penalisation. This is particularly useful in testing or
model choice problems but can also be effective in estimation problems, see for instance \citet{RSSB:RSSB781}. Here we require that the base measure puts exponentially small mass at the boundaries. This induce enough penalization to achieve consistency of the posterior distribution of $f(0)$ and $f(L)$. The following Theorem gives consistency of the posterior distribution of $f$ at every point on the support of $f_0$ including the boundaries. 

\begin{theo}\label{th2:allsupp}
	Let $x$ be in $[0,L]$ with with $0<L\leq \infty$ but $x < \infty$. Let $f_0 \in \F_L$ such that $f_0'$ exists at $x$ and $f_0'(x) < 0$. Let $X_i$ , $i=1, \dots ,n $ and $\Pi$ be
	either a Type I or Type II prior satisfying \textbf{C2} or \textbf{C3} with $\alpha$ satisfying condition \eqref{eq:cond:alpha:b} if $x=0$ or \eqref{eq:cond:alpha:c} if $x = L$. Then, for all $x$ in $[0,L]$ with $x< \infty$, and $\epsilon >0$  
\begin{equation}
{\Pi}\big(|f_P(x) - f_0(x)|> \epsilon | \Xn \big) \to 0  .
\end{equation}
Consider the posterior median $\fha_n^\pi(x) = \inf \{ t, \Pi\big[f_P(x) \leq t | \Xn\big] > 1/2 \}$, it follows that  
\begin{equation}
P_0\big(|\fha^\pi_n(x) - f_0(x)| > \epsilon | \Xn \big) \to 0  .
\end{equation} 
\end{theo}

The problem of estimating $f_0(0)$
under monotonicity constraints is another example of the effectiveness of penalisation induced by integration on the parameters. Although we do not have a proof for inconsistency of the posterior of $f(0)$ or $f(L)$ when $\alpha$ satisfies \eqref{eq:cond:alpha:a}, we believe that the similarly to the maximum likelihood estimator, the posterior distribution is in this case not consistent.

The following Theorem gives an upper bound on the concentration rate of the posterior distribution under the pointwise loss.
\begin{theo}
	Let $f_0$ be in $\F_L$ with $0<L \leq \infty$ and $\Pi$ be either a Type I or Type II prior satisfying \textbf{C1} or \textbf{C2} respectively with $\alpha$ satisfying \textbf{C1}, and let $x$ be in $(0,L)$ such that $f'$ exists in a neighbourhood of $x$ and $f'(x) < 0$, then for $C$ a positive constant
	\begin{equation}
		{\Pi}\left(|f_P(x) - f_0(x)|> C\left(n\over \log(n)  \right)^{-2/9} | \Xn \right) \to 0  .
	\end{equation}
when $n$ goes to infinity. 
\label{th:cvrate_fixedpoint}
\end{theo}

Here the concentration rate is subobtimal. It is however the best rate that one can obtain using the usual approach by testing \citep[see][]{Hebert2013concent} . Proving that the posterior concentrates at the rate $n^{-1/3}$ up to some power of $\log(n)$ would require some more refined control of the posterior distribution close to Bernstein von Mise types of results, see \citet{castillo2012semiparametric}, which in the case of mixture models is very difficult and beyond the scope of this paper. 

We derive from Theorem \ref{th2:allsupp} the consistency of the posterior distribution for the sup norm. This is particularly useful when considering confidence bands, as pointed out in \citet{GineNickl2010}. Under similar assumptions as in \citet{Durot2012}, we get the consistency of the posterior distribution for the sup norm. Note that contrariwise to \citet{Durot2012}, we do not restrict to sub-intervals of the support of the density. This is mainly due to the fact that the Bayesian approaches are consistent at the boundaries of the support of $f_0$. 

\begin{theo} 
	Let $f_0 \in \F_L$ with $0<L< \infty$ be such that $f_0'$ exists and $||f_0'||_\infty < \infty$ and for all $x\in [0,L]$, $f_0'(x) < 0$. Let also the prior $\Pi$ be either a Type I or Type II prior satisfying \textbf{C1} or \textbf{C2} with $\alpha$ satisfying conditions \eqref{eq:cond:alpha:b} and \eqref{eq:cond:alpha:c} respectively. Then 
\begin{equation}
\Pi( \sup_{x \in [0,L]}|f_P(x)-f_0(x)| > \epsilon | X_n) \to 0  .
\end{equation} 
\label{th3}
\end{theo}

Similar results as in Theorem \ref{th:cvrate_fixedpoint} also hold for the concentration rate of the posterior distribution for the supremum over all subsets of the form $(a,b)$ with $0<a<b<L$ with the same rate. 
\section{Proofs} \label{sec:proof}
In this section we prove Theorems \ref{th1} to \ref{th3} given in Section \ref{sec:main}. 
To prove Theorems 3-6, we need to construct tests that are adapted to the pointwise or supremum loss. The usual approach based on \cite{cam1986asymptotic} cannot be applied in this case. We thus construct test based on the Maximum Likelihood Estimator. 

\subsection{Proof of Theorems \ref{th1} and \ref{th15} }
The proofs of Theorems \ref{th1} and \ref{th15} follow the general ideas of \citet{ggv00} with some modification due to the fact that the Kullback-Leibler property is not satisfied.
We first focus on density on $\F_L$ with $L< \infty$ and extend these results to monotone non increasing density with support $\R^+$ that satisfy C3.
We extended the approach used in \citet{KHAZAEI:2010:INRIA-00494692:1} to the concentration rate framework and get similar results as those presented in \citet{ggv00}. More precisely, the proofs relies on the following Theorem which is a modification of \citet{ggv00} main Theorem proposed by \cite{RouRivBVM}. To tackle the fact that the usual Kullback Leibler property is not satisfied in its usual sense, we consider truncated versions of the densities 
\begin{equation} \label{eq:def_fn}
f_n(\cdot) = {f(\cdot) \I_{[0,\theta_n]}(\cdot) \over F(\theta_n) } ,\;  f_{0,n}(\cdot) = {f_0(\cdot) \I_{[0,\theta_n]}(\cdot) \over F_0(\theta_n) }
\end{equation}
where $\theta_n$ is defined as 
$$
\theta_n = \inf\{x, 1-F_0(x) < \frac{\epsilon_n}{2n}\}.
$$
We then define the counterpart of the Kullback Leibler neighbourhoods 
\begin{multline}
\label{eq:def_Sn}
S_n(\epsilon_n,\theta_n) =  
 \Bigg\{f, KL(f_n,f_{0,n}) \leq \epsilon_n^ 2, \\ \int f_{0,n}(x) \left (\log\left({f(x) \over f_0(x)} \right)\right )^2dx \leq \epsilon_n^2 ,  \int_0^{\theta_n} f(x) dx 
\gtrsim 1-  \epsilon_n^2 \Bigg\}.
\end{multline}
\begin{theo}\label{th:RivRous}
Let $f_0$ be the true density and let $\Pi$ be a prior on $\F$ satisfying the following conditions : there exist a sequence $(\epsilon_n)$ such that $\epsilon_n \to 0$ and $n\epsilon_n^2 \to \infty$ and a constant $c > 0$ such that for any $n$ there exist  $\F_n \subset \F$ satisfying 

$$ 
\Pi(\F_n^c) = o (\exp(-(c+2)n\epsilon_n^2)) .
$$

For any $j \in \N$, $j>0$, let $\F_{n,j} = \{ f \in \F_n, j\epsilon_n < d(f,f_0)  \leq (j+1) \epsilon_n \}$ and $N_{n,j}$ the Hellinger (or $L_1$) metric entropy of $\F_{n,j}$. There exists a $J_{0,n}$ such that for all $j \geq J_{0,n}$ 

$$
N_{n,j} \leq (K-1) n\epsilon_n^2 j^2 ,
$$

where $K$ is an absolute constant. \\ 

Let $S_n(\epsilon_n, \theta_n)$ be defined as in \eqref{eq:def_Sn} and let $\Pi$ be such that
\begin{equation}
\Pi(S_n(\epsilon_n,\theta_n)) \geq \exp(-cn\epsilon_n^2)   .
\label{eq:cond_kull}
\end{equation}
We have : 
$$
\Pi(f : d(f_0,f) \leq J_{0,n} \epsilon_n | \Xn) = 1 + {o}_P(1) .
$$
\end{theo}
The proof of this Theorem is postponed to Appendix \ref{app:B}. We will thus prove that the conditions of Theorem \ref{th:RivRous} are satisfied in our case. 
Let $f_0$ be in $\F_L$. The following lemma states that \eqref{eq:cond_kull} is satisfied.
\begin{lemm}\label{lem:piSn}
Let $\Pi$ be either a Type 1 or Type 2 prior on $\F_L$ as in Theorem \ref{th1} and let $S_n(\epsilon_n,\theta_n)$ be a set as in \eqref{eq:def_Sn}, then 
\begin{equation}
\Pi(S_n(\epsilon_n, \theta_n) )\gtrsim \exp\Bigg\{C_1\epsilon_n^{-1}\log(\epsilon_n)\Bigg\} .
\end{equation}
\end{lemm}
This lemma is proved in appendix \ref{app:A}.
The $\epsilon$ metric entropy of the set of bounded monotone non increasing densities has been shown to be less than $\epsilon^{-1}$, up to a constant (see \cite{MR873578} or \citet{WnVwcep} for instance). As the prior puts mass on $\F_L$, on which $f(0)$ is not uniformly bounded, we consider an increasing sequence of sieves
\begin{equation}
\F_n = \big\{ f\in \F_L, f(0) \leq M_n \Big\} .
\label{eq:Fn}
\end{equation}
 where $M_n = \exp\Big\{ cn^{1/3} \log(n)^{2/3}(t_2+1)^{-1} \Big\}$ with $t_2$ as in the conditions C1 or C2. The following Lemma shows that $\F_n$ covers most of the support of $\Pi$ as $n$ increase. 
\begin{lemm}
Let $\F_n$ be defined by \eqref{eq:Fn} and $\Pi$ be either a Type $1$ or Type $2$ as in Theorem \ref{th1}, then 
$$
\Pi\big(\F_n^c \big) \lesssim e^{-cn^{1/3} \log(n)^{2/3}} .
$$
\label{lem:Fn}
\end{lemm}
Here again, the proof is postponed to appendix \ref{app:A}.
We now get an upper bound for the $\epsilon$-metric entropy of the set $\F_n$. Recall that in \citet{MR822052} it is proved that the $L_1$ metric entropy of monotone non increasing densities on $[0,1]$ bounded by $M$ can be bounded from above
by $C_0\log(M) \epsilon_n^{-1}$. We cannot apply this result directly for the sets $\F_n$ as it would give a suboptimal control of the entropy to construct tests in a similar way as in \citet{ggv00}. In fact the upper bound on the entropy of $\F_n$ is of the order of $e^{n\epsilon_n}$ the usual conditions of \citet{ggv00} requires an upper bound of the order $e^{n\epsilon_n^2}$. However as stated in Theorem \ref{th:RivRous} it is enough to bound the $\epsilon$-metric entropy of the sets 
$$
\F_{n,j} = \left \{ f \in \F_n, j\epsilon_n \leq d(f,f_0) \leq (j+1) \epsilon_n \right  \} ,
$$
for $j \in \N^*$. We can easily adapt the results of 
\citet{MR822052} to positive monotone non increasing functions on any interval $[a,b]$ and get the following Lemma. 

\begin{lemm}
Let $\tilde{\F}$ be the set of positive monotone non increasing functions on $[a,b]$ such that for all $f$ in $\tilde{\F}$,$\int_a^b f \leq M_2 $ and $f \leq M$, then
$$
N(\epsilon,\tilde{\F},d) \lesssim \epsilon^{-1} \log(M+1) \Big((b-a) + 3M_2 \Big) .
$$ 
\label{lem:entropGroe}
\end{lemm}
The proof of this Lemma is straightforward given the results of \cite{MR822052} and is thus omitted. 
Let $x_{n,j} \in [0,L]$ such that $\epsilon_n/2 \leq x_{n,j} \leq \epsilon_n$. We denote for all $f$ in $\F_{n,j}$ $f_{1,j} = f \I_{[0,x_{n,j})}$ and $f_{2,j} = f \I_{[x_{n,j},L]}$. Since for all $f$ in $\F_{n,j}$ we have $\int_0^1 |f(x) - f_0(x) | dx \leq (j+1) \epsilon_n $ then 

$$ 
\int_0^{x_{n,j}} f(x) dx - \int_0^{x_{n,j}} f_0(x) dx \leq  (j+1) \epsilon_n  ,
$$
which implies that 
$$
x_{n,j}f(x_{n,j}) \leq  x_{n,j}f_0(0) + (j+1) \epsilon_n ,
$$
which in turn gives 
$$ 
f(x_{n,j}) \leq f_0(0) + 2 (j+1).
$$ 
 
Recall that for all $f \in \F_n$ we have $f(0) \leq M_n$. 
Using Lemma \ref{lem:entropGroe}, we construct an $\epsilon_n/2$-net for the set $\F_{n,j}^1 = \Big\{ f_{1,j}, f\in \F_{n,j} \Big\}$ with $N_1$ points, and 
$$
\log(N_1) \lesssim \epsilon_n^{-1} \log(M_n +1)   \epsilon_n (j+2),
$$
and thus deduce
\begin{equation}\label{eq:N1}
\log(N_1) \leq C' n \epsilon_n^2 j^2
\end{equation}  
Similarly, given that $f(x_{n,j}) \leq M + 2(j+1)$ we get an $\epsilon_n/2$-net for the set $\F_{n,j}^2 = \Big\{ f_{2,j}, f\in \F_{n,j} \Big\}$ with $N_2$ points and 

\begin{equation}\label{eq:N2}
\log(N_2) \leq \tilde{C}' n\epsilon_n^2 j^2 .
\end{equation}

This provide a $\epsilon_n$-net for $\F_{n,j}$ with less than $N_1\times N_2$ points. Given \eqref{eq:N1} and \eqref{eq:N2} the $L_1$ metric entropy of the sets $\F_{n,j}$ satisfy 

\begin{equation}\label{eq:entropy}  
\log(N(\F_{n,j},\epsilon_n,L_1)) \lesssim n\epsilon_n^2 j^2 .
\end{equation}

The conditions of Theorem \ref{th:RivRous} are thus satisfied which ends the proof of Theorem \ref{th1}



\paragraph*{Extention to $\R^+$}

Given that $ f_0(x) \lesssim e^{-\beta x^\tau}$ when $x$ goes to infinity, if $\theta_n$ is
such that $\theta_n = \inf \{ x, 1-F_0(x) < \epsilon_n/(2n) \}$ then $\theta_n \lesssim (\log(n))^{1/\tau}$. Using similar arguments as before, Lemma \ref{lem:piSn} still holds under the exponential tail assumption. We now get an upper bound for the $\epsilon$-metric entropy of $\F_{n,j}$. Here again, we split $\F_{n,j}$ into two parts. The construction of an $\epsilon_n/2$-net for $\F_{n,j}^1$ does not change and therefore \eqref{eq:N1} holds. Finally, let $\tilde{\F}_{n,j}^2 = \{ f \in \F_{n,j}^2, \forall x>\theta_n, f(x) = 0 \}$. Given Lemma \ref{lem:entropGroe}, we get for $c_1>0$ large enough an $\epsilon_n/(2c_1(j+1))$-net for $\tilde{\F}_{n,j}^2$ by considering $f^\star$ the restriction of $f$ to $[x_{n,j}, \theta_n]$. We have 
$$
d(f,f^\star) \leq c_2 (j+1) \epsilon_n ,
$$
where $d(\cdot,\cdot)$ is either the $L_1$ or Hellinger distance. Hence, for $c_1 > c_2$ an $\epsilon/2$-net for $\F_{n,j}^2$ with at most $e^{c_3n\epsilon_n^2j^2}$ points and thus 
$$
\log\left ( N(\F_{n,j}^2,\epsilon_n,d \right ) )\leq \tilde{C}'' n \epsilon_n^2j^2.
$$ 
We conclude using the same arguments as in the preceding section, and thus Theorem \ref{th15} is proved.

\subsection{Proof of Theorems \ref{th2} and \ref{th:cvrate_fixedpoint} } 
\label{sec:proof:2}

To prove Theorem \ref{th2} and \ref{th:cvrate_fixedpoint}, we need to construct tests for all $x \in (0,L)$ of $f_0$ versus $|f_P(x) - f_0(x)| \geq \epsilon_n^{2/3}$ as the approach used in \citet{ggv00} is not suited for the pointwise loss. As we have $\Pi(||f_P - f_0||_1 > \epsilon_n |\Xn) = o_{P_0}(1)$ we can consider functions $f_P$ such that $||f_P-f_0||_1 \leq \epsilon_n$. We construct tests $\Phi_n$ such that 
$$
E_0^n(\Phi) = o(1) ,~ \sup_{f, |f(x) - f_0(x)|> \epsilon_n} \E_f^n(1-\Phi) \leq e^{-Cn\epsilon_n^2}.
$$

Denote $A_\epsilon^x := \{ f, |f(x) - f_0(x)| > \epsilon\}$ that can be split into $A_\epsilon^{x,+} = \{ f, f(x) - f_0(x) > \epsilon\}$ and $A_\epsilon^{x,-} = \{ f, f(x) - f_0(x) < - \epsilon\}$ and denote $e_n = e_0 \epsilon_n^{2/3}$ and $h_n = h_0 e_n$. Consider the tests 

\begin{eqnarray*}
\phi^+_n &=& \I\left\{ n^{-1} \sum_{i=1}^n \I_{[x-h_n,x]} (X_i) - \int_{x-h_n}^x f_0(t) dt > c_n\right\} \\ 
\phi^-_n &=& \I\left\{ n^{-1} \sum_{i=1}^n \I_{[x,x+h_n]} (X_i) - \int_{x}^{x+h_n} f_0(t) dt < - c_n\right\}
\end{eqnarray*}

We immediately get $\E_0^n(\max(\phi_n^+, \phi_n^-) = o(1)$.
Note that if $f_P(x) > f_0(x) + e_n$ then 

\begin{align*}
\int_{x-h_n}^x f_P(t) - f_0(t) dt &\geq h_n(f_P(x) - f_0(x)) - \int_{x-h}^x f_0(t) - f_0(x)dt \\ 
&\geq h_ne_n -C_0h^2
\end{align*}

for some $C_0>0$ that only depends on $f_0$. Similarly if $f_P(x) < f_0(x) - e_n$ then for all $h>0$ 
$$
\int_{x}^{x+h} f_P(t) - f_0(t) dt \leq -he_n + C_0 h^2
$$

We thus deduce for $f_P$ such that $f_P(x) -f_0(x) > e_n$ 
\begin{align*}
P_f(1-\phi^+_n) &\leq P_f \left( n^{-1} \sum_{i=1}^n \I_{[x-h_n,x]}(X_i) - \int_{x-h_n}^x f_P(t) dt \leq -h_ne_n + C_0h^2 + c_n \right) \\ 
&\leq P_f\left( n^{-1} \sum_{i=1}^n \I_{[x-h,x]}(X_i) - \int_{x-h}^x f_P(t) dt \leq -h_0e^2_n / 2 \right),\\ 
\end{align*}
if $c_n\leq e_n^2$ and $h_0 \leq 1/C_0$. Now note that for $f_P$ such that $||f_P - f_0||_1 \leq \epsilon_n$  
\begin{align*}
\int_{x-h_n}^x f_P &\geq -\int_{0}^\infty |f - f_0| + \int_{x-h_n}^x f_0 \\
   &\geq  -\epsilon_n + \int_{x-h_n}^x f_0 \\ 
    &\geq -e_n + h_nf_0(x) \geq h_nf_0(x)/2 .
\end{align*}
Moreover, 

$$
\int_{x-h_n}^x f_P \leq e_n + h_nf_0(x-h_n) \leq 2h_nf_0(x)
$$
for $n$ large enough and $h$ small anough. We conclude that 
$$
\mathrm{Var}_{f_P}^n \left( n^{-1} \sum_{i=1}^n \I_{[x-h,x]}(X_i) \right) \leq 2 h f_0(x)
$$
Thus using Bernstein's inequality (e.g. \cite{WnVwcep} Lemma 2.2.9 p. 102) we get 

$$
P_f(1-\phi^+) \leq 2 e^{-nh_ne_n^2/(2+e_n/3)}.
$$

Similarly, we have 
$$
P_f(1-\phi^-_n) \leq 2 e^{-nh_ne_n^2/(2+e_n/3)}. 
$$

Taking $\Phi_n = \max(\phi^+_n, \phi^-_n)$ we deduce

\begin{align*}
P_0(\Phi_n) &= o(1) \\ 
\sup_{f\in A^x_{e_n}} P_f(1-\Phi_n) &\leq e^{-Ch_0e_n^3}
\end{align*}

We have 

\begin{align*}
P_0(\Phi_n) &= o(1) \\ 
\sup_{f\in A^x_{e_n}} P_f(1-\Phi_n) &\leq e^{-Cne_0\epsilon_n^2}
\end{align*}

Similarly to the proof of Theorem \ref{th:RivRous}, following \cite{KHAZAEI:2010:INRIA-00494692:1}, we get an
exponentially small lower bound for $D_n$. More precisely, we get that

$$
D_n \geq 2 e^{-(c+2)n \epsilon_n^2}
$$
with probability that goes to $1$. 
Note that
\begin{equation}
\begin{split}
\mathbf{E}_0^n\left( {N_n \over D_n} \right) &\leq \mathbf{E}_0^n (\Phi_n^x) +P_0^n(D_n \leq  e^{-(c+2)n \epsilon_n^2}) +
\\& \mathbf{E}_0^n(\Pi[\F_n^c | \Xn]) + e^{(c+2)n \epsilon_n^2}  \int_{A_\epsilon \cap \F_n}\mathbf{E}_f^n(1- \Phi_n^x) d\Pi(f)
\end{split}
\label{eq:ENnDn}.
\end{equation}

Given the preceding results, we have 

$$
\mathbf{E}_0^n\left( {N_n \over D_n} \right) \leq o(1) + e^{(c+2)n \epsilon_n^2}  \sup_f \mathbf{E}_f^n(1- \Phi_n^x)
$$

which ends the proof choosing $e_0$ large enough. 
\paragraph{Consistency of a Bayesian estimator}

We consider in this section $\fha^\pi_n(t)$, the Bayesian estimator associated with the absolute error loss, define as the median of the posterior distribution. Consistency of the posterior mean, which is the most common Bayesian estimator is however not proved here but could nevertheless be an interesting result.  \\ 

We first define $\fha^\pi_n(t)$ such that 

\begin{equation}
\fha_n^\pi(t) = \inf \{ x, \Pi [f_P(t) \leq x | \Xn ] > 1/2 \}.
\end{equation}

In order to get consistency in probability we note that if $\fha^\pi_n(t) - f_0(t) > \epsilon$ then 
$$
\Pi( f_P(t) > f_0(t) + \epsilon | \Xn) > 1/2.
$$ 

And if $\fha^\pi_n(t) - f_0(t) < - \epsilon$ then 
$$
\Pi( f_P(t) < f_0(t) - \epsilon | \Xn) > 1/2.
$$

We deduce, with Markov inequality and Theorem~\ref{th2}

\begin{eqnarray*}
P_0^n(\fha^\pi_n(t) - f_0(t)>\epsilon) &\leq& P^n_0(\Pi(f_P(t) > f_0(t) + \epsilon | \Xn) > 1/2 )\\
						   &\leq& 2\mathbf{E}^n_0(\Pi(f_P(t) > f_0(t) + \epsilon | \Xn) > 1/2 )\\ 
						   &\leq& o(1),
\end{eqnarray*}

and similarly 

$$ 
P_0^n(\fha^\pi_n(t) - f_0(t) < -\epsilon) \leq o(1).
$$

Thus we have $P_0^n(|\fha^\pi_n(t) - f_0(t)| > \epsilon) \to 0 $ which gives the consistency in probability of $\fha^\pi_n(t)$.
\subsection{Proof of Theorem \ref{th2:allsupp} }
The previous proof holds for all $x \in (0,L)$ we now need to prove the consistency of the posterior for $x=0$ and $x = L$, when the prior satisfies conditions \eqref{eq:cond:alpha:b} or \eqref{eq:cond:alpha:c}. We first consider the case $x=0$, the case $x=L$ can be deduce with symmetric arguments. 

As before, consider the set $A_\epsilon^0$ and split it in $A_\epsilon^{0,+}$ and $A_\epsilon^{0,-}$. Note that using the same test $\phi^-_n$ as before we easily get 
$$
\Pi( A_\epsilon^{0,-} | \Xn) = o_{P_0}(1).
$$
We now consider $f_P \in A_\epsilon^{0,+}$. As before we can restrict ourselves to functions $f_P$ such that $|| f_P - f_0 ||_1 \leq \epsilon_n$. We thus have for $h = 2\epsilon_n/\epsilon$ 
\begin{align*}
f_P(0) - f_0(0) &\leq f_P(0) - f_P(h) + h^{-1} \int |f_0(t) - f_P(t)| dt \\ 
& \leq f_P(0) - f_P(h) + h^{-1} \epsilon_n  \\
& = f_P(0) - f_P(h) + \epsilon/2 .
\end{align*}
We now prove that the prior mass of the event $\{ f_P(0) - f_P(h) > \epsilon/2\}$ is less that $e^{-(c+2)n\epsilon_n^2}$. Using Markov inequality we get 
$$
\Pi(f_P(0) - f_P(h) > \epsilon/2) \leq 2 \epsilon^{-1} \int_0^h \frac{1}{\theta} \alpha(\theta) d\theta \leq e^{-a_2/h} \lesssim e^{-a_2 n\epsilon_n^2\log(n)}.
$$
Using the same control for $D_n$ as in the proof of Theorem \ref{th:RivRous}, and applying the usual method of \citet{ggv00}, we get the desired result.
\subsection{Proof of Theorem \ref{th3} }
In this section we prove that the posterior distribution is consistent in sup norm.
Here again, the main difficulty is to construct tests that are adapted to the considered loss. More precisely we construct a test $\Phi$ such that 
$$
E_0^n(\Phi) = o(1) ,~ \sup_{f, \sup_{[0,L]}|f - f_0|> \epsilon_n} \E_f^n(1-\Phi) \leq e^{-Cn\epsilon_n^2}.
$$

 To do so we consider a combination of the tests considered in the previous section noting that if the posterior distribution is consistent at the points of a sufficiently refined partition of $[0,L]$ then it is consistent for the sup norm. Here again, we will only consider the case
$L=1$ without loss of generality. 
We first  denote 
$$
B_\epsilon = \left\{ f, \sup_{[0,L]}\{|f(x) - f_0(x)| > \epsilon \right\} 
$$
Let $C'_0$ be a positive constant such that $||f'_0||_\infty \leq C'_0$ and let $(x_i)_{i}$ be the separation points of a $\epsilon/(8C_0')$ regular partition of $[0,1]$ and $p = \mathrm{Card}\{ (x_i)_i\}$. Note that 
$$ 
B_\epsilon = \bigcup_{i=1}^p \{f, \sup_{[x_i,x_{i+1}]}\{|f(x) - f_0(x)| > \epsilon\}.
$$ 
Recall that $A_\epsilon^x = \{ f, |f(x) - f_0(x)| > \epsilon\}$. We consider the set
$
B_\epsilon \bigcap_{i=1}^p (A_{\epsilon/8}^{x_i}) ^c 
$. Given Theorem \ref{th2}, we have that 
$$
\E_0^n\left (\Pi\left (\bigcup_{i=1}^p (A_{\epsilon/5}^{x_i}) \Big|\Xn\right ) \right ) = o(1) .
$$
If $f \in B_\epsilon$ we have for all $x \in [x_i, x_{i+1}]$, 
$$
|f(x) - f_0(x)| \leq |f(x) - f(x_i)| + |f(x_i) - f_0(x_i)| + |f_0(x_i)- f_0(x)|. 
$$
Given that $f$ is monotone non increasing, and given the hypotheses on $f_0$ we have 
\begin{align*}
|f(x) - f(x_i)| &\leq |f(x_{i+1}) - f(x_i)| \\
& \leq |f(x_{i+1}) - f_0(x_{i+1})| + |f_0(x_{i+1}) - f_0(x_i)| + |f_0(x_i)- f(x_i)| \\
& \leq 3\epsilon/5
\end{align*}
and for the same reasons 
$$
|f(x_i) - f_0(x_i)| + |f_0(x_i)- f_0(x)| \leq 2\epsilon/5.
$$
Which leads to 
$$
|f(x) - f_0(x)| \leq \epsilon 
$$
and thus, taking the supremum over $x$, we get 
$$
\sup_{x \in [x_i,x_{i+1} ]} |f(x) - f_0(x)| \leq \epsilon. 
$$
We then deduce
$$
\Pi(B_\epsilon | \Xn) \leq \Pi\left (B_\epsilon \bigcap \left \{ \bigcap_{i=1}^p (A_{\epsilon/5}^{x_i})^c \right \} \right )  + \Pi\left (\bigcup_{i=1}^p (A_{\epsilon/5}^{x_i}) \right ) = o_{P_0}(1)
$$
Which gives the consistency of the posterior distribution in sup norm \\

\section{Discussion}

In this paper, we obtain an upper bound for the concentration rate of the posterior distribution under monotonicity constraints. This is of interest as in this model, the standard approach based on the seminal paper of \citet{ggv00} cannot be applied directly. We prove that the concentration rate of the posterior is (up to a $\log(n)$ factor) the minimax estimation rate $(n/\log(n))^{-1/3}$ for standard losses such as $L_1$ or Hellinger. 

We also prove that the posterior distribution is consistent for the pointwise loss at any point of the support and for the sup norm loss. Studying asymptotic properties for these losses is difficult in general as the usual approach are well suited for losses that are related to the Hellinger metric. Obtaining more refined results on the asymptotic behaviour of the posterior distribution will require refined control of the likelihood which in the case of nonparametric mixture models is a difficult task.  


\bibliographystyle{apalike} 
\bibliography{vitesses_final_mref}

\begin{thebibliography}{}

\bibitem[Balabdaoui and Wellner, 2007]{balabdaoui:wellner}
Balabdaoui, F. and Wellner, J.~A. (2007).
\newblock {Estimation of a {$k$}-monotone density: limit distribution theory
  and the spline connection}.
\newblock {\em Ann. Statist.}, 35(6):2536--2564.

\bibitem[Brunner and Lo, 1989]{BrLo1989}
Brunner, L.~J. and Lo, A.~Y. (1989).
\newblock {Bayes methods for a symmetric unimodal density and its mode}.
\newblock {\em Ann. Statist.}, 17(4):1550--1566.

\bibitem[Castillo, 2013]{castillo2012semiparametric}
Castillo, I. (2013).
\newblock On bayesian supremum norm contraction rates.
\newblock {\em arXiv preprint arXiv:1304.1761}.

\bibitem[Donoho and Liu, 1991]{donoho1991geometrizing}
Donoho, D.~L. and Liu, R.~C. (1991).
\newblock Geometrizing rates of convergence, ii.
\newblock {\em The Annals of Statistics}, pages 633--667.

\bibitem[Durot et~al., 2012]{Durot2012}
Durot, C., Kulikov, V.~N., Lopuha{\"a}, H.~P., et~al. (2012).
\newblock The limit distribution of the $l_infty$-error of grenander-type
  estimators.
\newblock {\em The Annals of Statistics}, 40(3):1578--1608.

\bibitem[Ferguson, 1983]{MR736538}
Ferguson, T.~S. (1983).
\newblock {Bayesian density estimation by mixtures of normal distributions}.
\newblock In {\em {Recent advances in statistics}}, pages 287--302. Academic
  Press, New York.

\bibitem[Ghosal et~al., 2000]{ggv00}
Ghosal, S., Ghosh, J.~K., and van~der Vaart, A.~W. (2000).
\newblock {Convergence rates of posterior distributions}.
\newblock {\em Ann. Statist.}, 28(2):500--531.

\bibitem[Ghosal and van~der Vaart, 2007]{GVdV07:Dirichlet}
Ghosal, S. and van~der Vaart, A. (2007).
\newblock {Posterior convergence rates of {D}irichlet mixtures at smooth
  densities}.
\newblock {\em Ann. Statist.}, 35(2):697--723.

\bibitem[Ghosal and van~der Vaart, 2001]{vdv2001}
Ghosal, S. and van~der Vaart, A.~W. (2001).
\newblock {Entropies and rates of convergence for maximum likelihood and
  {B}ayes estimation for mixtures of normal densities}.
\newblock {\em Ann. Statist.}, 29(5):1233--1263.

\bibitem[Gin{\'e} and Nickl, 2010]{GineNickl2010}
Gin{\'e}, E. and Nickl, R. (2010).
\newblock {Confidence bands in density estimation}.
\newblock {\em Ann. Statist.}, 38(2):1122--1170.

\bibitem[Grenander, 1956]{Gre1956}
Grenander, U. (1956).
\newblock {On the theory of mortality measurement. {II}}.
\newblock {\em Skand. Aktuarietidskr.}, 39:125--153 (1957).

\bibitem[Groeneboom, 1985]{MR822052}
Groeneboom, P. (1985).
\newblock {Estimating a monotone density}.
\newblock In {\em {Proceedings of the {B}erkeley conference in honor of {J}erzy
  {N}eyman and {J}ack {K}iefer, {V}ol.\ {II} ({B}erkeley, {C}alif., 1983)}},
  {Wadsworth Statist./Probab. Ser.}, pages 539--555, Belmont, CA. Wadsworth.

\bibitem[Groeneboom, 1986]{MR873578}
Groeneboom, P. (1986).
\newblock {Some current developments in density estimation}.
\newblock In {\em {Mathematics and computer science ({A}msterdam, 1983)}},
  volume~1 of {\em {CWI Monogr.}}, pages 163--192. North-Holland, Amsterdam.

\bibitem[Hoffmann et~al., 2013]{Hebert2013concent}
Hoffmann, M., Rousseau, J., and Schmidt-Hieber, J. (2013).
\newblock On adaptive posterior concentration rates.
\newblock {\em arXiv preprint arXiv:1305.5270}.

\bibitem[Khazaei et~al., 2010]{KHAZAEI:2010:INRIA-00494692:1}
Khazaei, S., Rousseau, J., and Balabdaoui, F. (2010).
\newblock {Bayesian Nonparametric Inference of decreasing densities}.
\newblock In {\em {42{\`e}mes Journ{\'e}es de Statistique}}, Marseille, France
  France.

\bibitem[Kruijer et~al., 2009]{kruijer:rousseau:vdv:09}
Kruijer, W., Rousseau, J., and van~der Vaart, A. (2009).
\newblock {Adaptive {Bayes}ian Density Estimation with Location-Scale
  Mixtures}.
\newblock Technical report.

\bibitem[{Le Cam}, 1986]{cam1986asymptotic}
{Le Cam}, L. (1986).
\newblock {\em {Asymptotic methods in statistical decision theory}}.
\newblock {Springer Series in Statistics}. Springer-Verlag, New York.

\bibitem[Lo, 1984]{lo84}
Lo, A.~Y. (1984).
\newblock {On a class of {B}ayesian nonparametric estimates. {I}. {D}ensity
  estimates}.
\newblock {\em Ann. Statist.}, 12(1):351--357.

\bibitem[{Prakasa Rao}, 1970]{PrakasaRao70}
{Prakasa Rao}, B. L.~S. (1970).
\newblock {Estimation for distributions with monotone failure rate}.
\newblock {\em Ann. Math. Statist.}, 41:507--519.

\bibitem[Reynaud-Bouret et~al., 2011]{RBRTM}
Reynaud-Bouret, P., Rivoirard, V., and Tuleau-Malot, C. (2011).
\newblock {Adaptive density estimation: a curse of support?}
\newblock {\em J. Statist. Plann. Inference}, 141(1):115--139.

\bibitem[Rivoirard et~al., 2012]{RouRivBVM}
Rivoirard, V., Rousseau, J., et~al. (2012).
\newblock Bernstein--von mises theorem for linear functionals of the density.
\newblock {\em The Annals of Statistics}, 40(3):1489--1523.

\bibitem[Robertson et~al., 1988]{robertson:wright:dykstra:1988}
Robertson, T., Wright, F.~T., and Dykstra, R.~L. (1988).
\newblock {\em {Order restricted statistical inference}}.
\newblock {Wiley Series in Probability and Mathematical Statistics: Probability
  and Mathematical Statistics}. John Wiley \& Sons Ltd., Chichester.

\bibitem[Rousseau, 2010]{Rousseau10}
Rousseau, J. (2010).
\newblock {Rates of convergence for the posterior distributions of mixtures of
  betas and adaptive nonparametric estimation of the density}.
\newblock {\em Ann. Statist.}, 38(1):146--180.

\bibitem[Rousseau and Mengersen, 2011]{RSSB:RSSB781}
Rousseau, J. and Mengersen, K. (2011).
\newblock {Asymptotic behaviour of the posterior distribution in overfitted
  mixture models}.
\newblock {\em J. R. Stat. Soc. Ser. B Stat. Methodol.}, 73(5):689--710.

\bibitem[Shen et~al., 2013]{SGmulti12}
Shen, W., Tokdar, S.~T., and Ghosal, S. (2013).
\newblock Adaptive bayesian multivariate density estimation with dirichlet
  mixtures.
\newblock {\em Biometrika}, 100(3):623--640.

\bibitem[Sun and Woodroofe, 1996]{MR1392133}
Sun, J. and Woodroofe, M. (1996).
\newblock {Adaptive smoothing for a penalized {NPMLE} of a non-increasing
  density}.
\newblock {\em J. Statist. Plann. Inference}, 52(2):143--159.

\bibitem[van~der Vaart and Wellner, 1996]{WnVwcep}
van~der Vaart, A.~W. and Wellner, J.~A. (1996).
\newblock {\em {Weak convergence and empirical processes}}.
\newblock {Springer Series in Statistics}. Springer-Verlag, New York.
\newblock With applications to statistics.

\bibitem[Williamson, 1956]{will:56}
Williamson, R.~E. (1956).
\newblock {Multiply monotone functions and their {L}aplace transforms}.
\newblock {\em Duke Math. J.}, 23:189--207.

\bibitem[Wu and Ghosal, 2008]{WG08}
Wu, Y. and Ghosal, S. (2008).
\newblock {Kullback {L}eibler property of kernel mixture priors in {B}ayesian
  density estimation}.
\newblock {\em Electron. J. Stat.}, 2:298--331.

\end{thebibliography}

\appendix 

\section{Technical Lemmas} 
\label{app:A}

\subsection{Proof of Lemma \ref{lem:piSn} }

To prove Lemma \ref{lem:piSn}, we first construct stepwise constant functions such that these approximations are in the truncated Kullback Leibler neighbourhood of $f_0$. We then construct a set $\mathcal{N}$ included in $S_n(\epsilon_n,\theta_n)$ based on the considered piecewise constant approximation such that for $\Pi$ a Type I or Type II prior $\Pi(\mathcal{N}) \geq e^{-Cn\epsilon_n^2}$.

We first construct a  piecewise constant approximation of $f_0$ which is base on a sequential subdivision of the interval $[0,L]$ with more refined subdivisions where $f_0$ is less regular such that the number of points is less than $\epsilon_n^{-1}$ points. 

This approximation is adapted from the proof of Theorem 2.5.7 in \citet{WnVwcep}.
We then identify a finite piecewise constant density by a mixture of uniform for which the Hellinger distance between the piecewise constant approximation $f_P$ of $f_0 \in \F$ and $f_0$ is less that $\epsilon_n$ and $||f_0/f_P||_\infty \leq M$.The following Lemma gives the form of a finite probability distribution $P$ such that $f_{P}$ is in the Kullback-Leibler neighbourhood of some $f\in \mathcal{F}$. 
\begin{lemm} \label{lem:KL}
Let $f \in \F_L$ be such that $f(0) \leq M < + \infty$. For all $0<\epsilon<1$ there exists $m \lesssim L^{1/3}M^{1/3} \epsilon^{-1}$, $p = (p_1, \dots, p_m) \in \mathfrak{S}_m$ and $x = (x_1,\dots, x_m) \in [0,L]^m$ such that $P = \sum_{i=1}^m \delta_{x_i} p_i$ satisfies 
\begin{equation}
KL(f,f_P) \lesssim \epsilon^2, \; \int  \left (\log\left({f\over f_P}  \right)\right )^2f \lesssim \epsilon^2,
\end{equation}
where $f_P$ is defined as in \eqref{eq:Will}.
\end{lemm}

\begin{proof}
For a fixed $\epsilon$, let $f$ be in $\F_L$. Consider $\P_0$ the coarsest partition : 
$$
0=x_0^0 < x_1^0 = L,
$$
at the $i^{th}$ step, let $\P_i$ be the partition 
$$ 
0=x_0^i<x_1^i<\cdots<x_{n_i}^i = L,
$$ 
and define 
$$ 
\eps_i = \max_j \left\{ (f(x_{j-1}^i) - f(x_j^i))(x_j^i - x_{j-1}^i)^{1/2}  \right\}.
$$
For each $j\geq 1$, if $
(f(x_{j-1}^i) - f(x_j^i))(x_j^i - x_{j-1}^i)^{1/2} \geq {\eps_i \over \sqrt{2}} $ we split the interval $[x_{j-1}, x_j]$ into two subsets of equal length. We then get a new partition $\P_{i+1}$. We continue the partitioning until the first $k$ such that $\eps_k^2 \leq \epsilon^3$.
At each step $i$, let $n_i$ be the number of intervals in $\P_i$, $s_i$ the number of interval in $\P_i$ that have been divided to obtain $\P_{i+1}$, and $c = 1/\sqrt{2}$. Thus, it is clear that $\eps_{i+1} \leq c\eps_i$ 

\begin{eqnarray*}
 s_i(c\eps_i)^{2/3} & \leq & \sum_j (f(x_{j-1}^i) - f(x_{j}^i) )^{2/3}(x_{j}^i-x_{j-1}^i) ^ {1/3} \\
                    & \leq & \left( \sum_j f(x_{j-1}^i) - f(x_{j}^i) \right)^{2/3}  \left( \sum_j x_{j}^i-x_{j-1}^i\right)^{1/3} 
                     \leq  M^{2/3}L^{1/3},
\end{eqnarray*}

using H\"older inequality. We then deduce that 

\begin{align*}
\sum_{j=1}^k n_j  =  k + \sum_{j = 1}^k j  s_{k-j}   \leq  2 \sum_{j = 1}^k j  s_{k-j} 
                & \leq 2 \sum_{j = 1}^k j  M^{2/3} L^{1/3} (c\eps_{k-j})^{- 2/3} \\
                & \leq  2  M^{2/3} L^{1/3} \eps_k^{-2/3}  2^{1/3} \sum_{j=1}^k j 2^{-j/3} \\
                & \leq   K_0  M^{2/3} L^{1/3} \eps_k^{-2/3},
\end{align*}
where $K_0 = 2(1-2^{-2/3})^{-2}$. Thus 
\begin{equation} \label{eq:nk}
n_k \leq K_0  M^{2/3} L^{1/3} \epsilon^{-1}.
\end{equation}

Now, for $f \in \F_L$, we prove that there exists a stepwise density with less than $K_0 M^{2/3}L^{1/3} \frac{1}{\epsilon}$ pieces such that 
\begin{equation} \label{eq:KL_fP}
KL(f,h) \leq \epsilon^2 \text{  and  } \int f \log({f_0 \over f_P})^2(x) dx \lesssim \epsilon^2
\end{equation}

In order to simplify notations, we define 

$$
x_i = x_i^k, \quad
l_i = x_i - x_{i-1}, \quad  
g_i = f(x_{i-1})^{1/2}.
$$

We consider the partition constructed above associated with $f^{1/2}$, which is also a monotone nonincreasing function that satisfy $f^{1/2}(0) \leq M^{1/2}$ (instead of $M$). We denote $g$ the function defined as $g(x) = \sum \I_{[x_{i-
1},x_i]}(x) g_i $  

\begin{equation*}
\begin{split}
||f^{1/2}-g ||_2^2 = \int (f^{1/2}-g)^2(x) dx &= \sum_{i=1}^{n_k} \int_{I_i} (f^{1/2}-g)^2(x) dx \\
&\leq  \sum_{i=1}^{n_k} \int_{I_i} (f^{1/2}(x_{i-1}^k) - f^{1/2}(x_{i}^k))^2 dx \\
&\leq  \sum_{i=1}^{n_k} (x_i^k - x_{i-1}^k) (f^{1/2}(x_{i-1}^k) - f^{1/2}(x_{i}^k))^2 \\
&\leq  n_k \eps_k^2 \leq  L^{1/3}K_0 M^{1/3} \epsilon^2 .\\
\end{split}
\end{equation*} 

We then define $h = \frac{g^2}{\int g^2} $ and and get an equivalent of $\int g^2$. 

\begin{align*}
 \int g^2 dx &= \int (g^2 - f)(x) dx + 1 \\ 
	     &= \int (g-\sqrt{f})(g+\sqrt{f})(x) dx +1\\
	     &= 1 + \O(\epsilon),
\end{align*}

and deduce that $(\int g^2)^{1/2} = 1 + \O(\eps)$. Let $H$ be the Hellinger distance  

\begin{align*}
H(f,h) &= H\left( f,\frac{g^2}{\int g^2} \right) \\ 
	 &\leq H(f,g^2) + H(g^2,\frac{g^2}{\int g^2}) \\
	 &\leq L^{1/6}K_0 M^{1/6} \epsilon+ \left(\int (g - \frac{g}{(\int g^2)^{1/2}})^2(x) dx \right)^{1/2} \lesssim \epsilon .\\
\end{align*}

Since $|| f/h||_\infty = ||f/g^2||_\infty (\int g^2) \leq (\int g^2) $, together with the above bound on $H(f,h)$ and Lemma 8 from \cite{GVdV07:Dirichlet}, we obtain the required result.

Let $P$ be a probability distribution defined by 
$$ 
P = \sum_{i=1}^{n_k} p_i \delta(x_i^k) \quad p_i = (h_{i-1} - h_{i})x_i^k \quad p_{n_k} = h_{n_k}x_{n_k}^k = h_{n_k}L,
$$
thus $f_P = h$ and given the previous result, lemma \ref{lem:KL} is proved. 
\end{proof}
Given Lemma \ref{lem:KL}, we now prove Lemma \ref{lem:piSn}.
\begin{proof}[Proof of Lemma \ref{lem:piSn}] We first consider the case where $ \theta^{t_1} \lesssim \alpha(\theta) \lesssim \theta^{t_2}$ for small $\theta$.
For $\epsilon_n$ as in Theorem \ref{th1}, define $\theta_n$ as 
$$
\theta_n = \inf \{ x, 1-F_0(x) < \frac{\epsilon_n}{2n} \} .
$$
Note that $F_0$ is c\`adl\`ag, thus 
\begin{equation}
F_0(\theta_n) \geq 1 - \epsilon_n/(2n)~ \text{and}~\forall y < \theta_n 1-F_0(y) > \epsilon_n/(2n).
\label{eq:Fthetan}
\end{equation}.
Using lemma \ref{lem:KL} with $L = \theta_n$, we obtain that there exists a distribution $P=\sum_{i=1}^{n_k} \delta_{x_i} p_i$ such that
$$
KL(f_{0,n},f_{P}) \leq  \epsilon_n^2, \text{ and } \int f_{0,n} \log\left({f_{0,n}\over f_{P}} \right)^2 \lesssim \epsilon_n^2.
$$
Note that $f_P$ has support $[0, \theta_n]$ and is such that $f_P(\theta_n) > 0$. Now, set $m = n_k$ and consider $P'$ the mixing distribution associated with $\{m,  x'_1, \dots, x'_m, p'_1  \dots,  p'_m \}$ with $\sum_{i=1}^{m} p'_i = 1$. Define for $1 \leq i \leq m-1$ the set $U_i = [0\vee (x_i - \epsilon_n^3/M, x_i + \epsilon_n^3/M]$ and $U_m = (\theta_n, \theta_n + \epsilon_n(L-\theta_n) \wedge \epsilon_n^3/M]$. Construct $P'$ such that $x'_i \in U_i$ and $|P'(U_i) - p_i| \leq \epsilon^2 m^{-1}$. 
We get 
$$
\forall t \in [0, \theta_n ] \; f_P'(t) > \frac{p'_m}{x'_m}.
$$
Given that $x'_m \in U_m$, we get $x'_m \leq \theta_n + \epsilon_n(L-\theta_n) \wedge \epsilon_n^3/M \lesssim \theta_n$ for $n$ large enough. Note also that $p'_m \geq p_m - \epsilon_n^2m^{-1}$. Given the construction of Lemma \ref{lem:KL}, we deduce 
$$
p_m \geq \frac{f_0(x_{i-1})}{1+\O(\epsilon_n)} \gtrsim f_0(x_{i-1}),
$$
for $n$ large enough. Furthermore, given \eqref{eq:Fthetan}
$$
\forall z< \theta_n, \; f_0(z) (L-z) \geq \int_z^L f_0(t) dt \geq {\epsilon_n \over 2n} ,
$$
thus 
$$
\forall t \in [0, \theta_n ] \; f_P'(t) \gtrsim { \frac{\epsilon_n}{2n} - \epsilon_n^2 m^{-1} \over \theta_n} \gtrsim \frac{\epsilon_n}{n}  ,
$$
and deduce that $||f_0/f_{P'}||_\infty \lesssim \frac{n}{\epsilon_n}$ 
Lemma 8 from \cite{GVdV07:Dirichlet} gives us that 
\begin{align*}
\int_{0}^{\theta_n} f_0(x) \log\left( {f_0  \over f_{P'}}\right)(x) dx &\lesssim \left(\epsilon_n^2 + H^2(f_{P},f_{P'})\right)(1+|\log(\epsilon_n/n)|) \\
&\lesssim \left(\epsilon_n^2 +|f_{P} -f_{P'}|_1 \right)(1+|\log(\epsilon_n/n)|) .
\end{align*}
Given the mixture representation \eqref{eq:Will} of $f_0$ and $f_P$, we get 
\begin{align*}
&\left(\epsilon_n^2 +|f_{P} -f_{P'}|_1 \right)(1+\log(n)) \\ 
&\lesssim \Big(\epsilon_n^2 + \int_0^{\theta_n} \Big| \sum(\frac{p_i}{x_i} - \frac{p'_i}{x'_i}) \I_{x \leq x_i} + \sum \frac{p_i}{x_i} (\I_{x \leq x_i} - \I_{x \leq x'_i})\Big|dx \Big)(1+\log(n)) \\ 
& \lesssim \Big( \epsilon_n^2 +  \sum |\frac{x_i}{x'_i} -1 |p'_i + \sum |p'_i - p_i| + \sum \frac{p_i}{x_i} |x'_i - x_i| \Big)(1+|\log(n)|) \\ 
&\lesssim \epsilon_n^2(1+|\log(n)|).
\end{align*}
Generally speaking, denoting $U_0 = [0,1] \cap \left( \cup_{i=1}^m U_i\right)^c$ and $\mathcal{N} = \{ P', |P'(U_i) - p_i|\leq \epsilon_n^2 m^{-1}\}$ we obtain that for all $P' \in \mathcal{N}$ 

$$
\int_0^{\theta_n} f_0(x)\log\Big(\frac{f_0}{f_{P'}}\Big)(x)dx \lesssim \epsilon_n^2 (1+|\log(n)|),
$$
and similarly 
$$
\int_0^{\theta_n} f_0(x)\log\Big(\frac{f_0}{f_{P'}}\Big)^2(x)dx \lesssim \epsilon_n^2 (1+|\log(n)|)^2,
$$
for $\epsilon_n$ small enough. Note also that for all $P' \in \mathcal{N}$ and $n$ large enough, as before we get 
$$
\int_{\theta_n}^L f_{P'}(x)dx \lesssim \frac{\epsilon_n}{n}.
$$

We now derive a control on $k$, the number of steps until $\varepsilon_k \leq \epsilon_n^{3/2}$ in the construction of Lemma \ref{lem:KL}. At step $k-1$, we have $\varepsilon_{k-1} \geq \epsilon_n^{3/2}$. It is clear that for all $j$, $\varepsilon_{j} \leq 2^{-1/2} \varepsilon_{j-1}$, thus 
\begin{equation*}
\begin{split}
M^{1/2}L^{1/2} 2^{-(k-1)/2} \geq \eps_{k-1} \geq \epsilon_n^{3/2} \\
\log(M^{1/2}L^{1/2}) - (k-1)\frac{\log(2)}{2} \geq \frac{3}{2}\log(\epsilon_n). \\
\end{split}
\end{equation*}
Finally, we have 
\begin{equation}\label{eq:k}
k \leq \frac{2}{\log(2)}(\log(M^{1/2}L^{1/2}) - \frac{3}{2}\log (\epsilon_n))  + 1 .\\
\end{equation}

We can then get a lower bound for $\Pi[\mathcal{N}]$ and, given that for $\epsilon_n$ small enough and $n$ large enough, we have 
$$
\mathcal{N} \subset S_n(\epsilon_n,\theta_n),
$$ 
we can deduce a lower bound for $\Pi\Big(S_n(\epsilon_n,\theta_n)\Big)$. For the Type $1$ prior, we have similarly to \citet{ggv00} 

\begin{eqnarray*}
\Pi[\mathcal{N}] &=& \Pr(\mathcal{D}(A\alpha(U_0), \dots,A\alpha(U_{n_k})) \in [p_i \pm \epsilon_n^2/n_k]) \\
		 &\geq& \frac{\Gamma(A)}{\prod_i \Gamma(A\alpha(U_i))} \prod_j \int_{(p_i - \epsilon_n^2/n_k)\wedge 0}^{(p_i +  \epsilon_n^2/n_k)} x_j^{A\alpha(U_j) - 1 } dx_j .\\
\end{eqnarray*}

Given condition C1, we have

$$
\alpha(U_i) \geq  \int_{U_i} \alpha_0 \theta^{t_1} d\theta, \\
$$
thus 
$$
\alpha(U_i) \geq 2\epsilon_n^3  \alpha_0 {x_i}^{t_1}.
$$

%
for $n$ large enough and $ \epsilon$ sufficiently small we have as in Lemma 6.1 of \citet{ggv00} 
$$
\Pi(\mathcal{N}) \gtrsim \exp\left\{ C_1 n_k \log(\epsilon) \right\} .
$$
Note that given \eqref{eq:nk}, $n_k \lesssim \epsilon_n^{-1}$ which gives the desired result. For the Type $2$ prior, we write

$$
\mathcal{N}' = \left\{ P' = \sum_{j=1}^{n_k} p'_j \delta_{x'_j}, |p'_j - p_j| \leq \epsilon^2 / n_k , |x'_j - x_j| \leq \epsilon_n^3  \right\} \subset S_n(\epsilon_n,\theta_n),
$$

we then deduce a lower bound for $\Pi[ S_n(\epsilon_n,\theta_n)]$

\begin{eqnarray*}
\Pi[\mathcal{N}'] &\geq& Q(K = n_k) \prod_{j=1}^{n_k} n_k^{-n_k} c^{n_k} \int_{\max(0,p_i - \epsilon^2 / n_k)}^{p_i + \epsilon^2 / n_k} w_j^{a_j} dw_j \prod_{j=1}^{n_k}\alpha(U_i) \\ 
		  &\geq& \exp \left\{  -cn_k\log{n_k} + \sum \log(\alpha(U_i))  + n_k \log(c) - n_k \log(n_k) + \sum a_j \log(2\epsilon^2 / n_k)\right\}  \\		 
		  &\gtrsim& \exp\left\{ C'_1 \epsilon^{-1} \log(\epsilon) \right\}.\\	
\end{eqnarray*}

We now consider the case where $e^{-a_1/\theta} \leq \alpha(\theta) \leq e^{-a_2/\theta}$ if $\theta$ is close to $0$ and $\sup_{x \in [0,\delta]} |f'_0(x)|\leq C_0$. We have that for $n$ large enough and $C>0$, a constant depending on $f_0$, $f_0(0) - f_0(\epsilon_n) \leq C \epsilon_n$.  Following Lemma \ref{lem:KL}, we can construct a piecewise constant approximation of $f_0$ on $[\delta,L]$. On $[0,\delta]$, consider the regular partition with $\lfloor \epsilon_n^{-1}\rfloor$ points and the piecewise constant approximation of $f_0$ defined as before (i.e. $f_i = f_0(x_{i-1})$). Again, this approximation can be identified with a measure $P$. Given the assumptions on $f_0$ we immediately get that $KL(f_0,f_P) \lesssim \epsilon_n^2$. 

 Consider the same sets $\mathcal{N}$ as before, with the same partitions $U_1, \dots, U_n$. Using similar computations as in Lemma $6.1$ of \citet{ggv00} we get that 
$$
\Pi(\mathcal{N}) \geq \exp\left\{ C_1 (n_k + \epsilon_n^{-1})  \log(\epsilon_n) + \sum \log(\alpha(U_i)) \right\} 
$$

For the $U_i$ included in $[\delta,L]$ we have $\alpha(U_i) \gtrsim \epsilon_n^{3/2}$. For the $U_i$ included in $[0,\delta]$ we have $\alpha(U_i) \gtrsim \epsilon_n \exp\left\{ -a/(i\epsilon_n)\right\}$, which gives 
$$
\sum \alpha(U_i) \lesssim -\epsilon_n^{-1} \log(n)
$$
We end the proof using similar argument as before. 
\end{proof}

\subsection{Proof of Lemma \ref{lem:Fn} } 
The proof of Lemma \ref{lem:Fn} is straightforward and comes directly from C1 and C2. 
\begin{proof}
Recall that given \eqref{eq:Will}, $f(0) = \int_{[0,1]} {1 \over \theta} dP(\theta)$. Then

$$
\Pi\left[ \int_0^1 \frac{1}{\theta} dP(\theta) \geq M_n  \right] = \Pi\left[ \int_0^{2M_n^{-1}} \frac{1}{\theta} dP(\theta) + \int_{2M_n^{-1}}^1 \frac{1}{\theta} dP(\theta) \geq M_n  \right].
$$ 

Note that 

$$
\int_{2M_n^{-1}}^1 \frac{1}{\theta} dP(\theta) \leq M_n/2 \int_{2M_n^{-1}}^1 dP(\theta) \leq M_n/2.
$$

Thus the set $\{ P, \int_0^{2M_n^{-1}} \theta^{-1} dP(\theta) \geq M_n/2\}$ contains $  \F_n^c$ and 

\begin{eqnarray*}
\Pi [\F_n^c] &\leq& \Pi \left[ \int_0^{2M_n^{-1}} {1 \over \theta} dP(\theta) > M_n /2 \right] \\ 
	     & \leq& 2M_n^{-1} E\left[ \int_0^{2M_n^{-1}} {1 \over \theta} dP(\theta)  \right],
\end{eqnarray*}

using Markov inequality. Then for a Type $1$ prior when $n$ large enough 

\begin{eqnarray*}
\Pi [\F_n^c] & \leq& 2M_n^{-1} \int_0^{2M_n^{-1}} {1 \over \theta} \alpha(\theta) d\theta  \\ 
	     & \leq& 2M_n^{-1} \int_0^{2M_n^{-1}} \theta^{t_2-1} d\theta = {(2M_n^{-1})^{t_2+1} \over t_2} = C e^{-cn^{1/3}\log(n)^{2/3}}.
\end{eqnarray*}

For a Type $2$ prior, we have that 

\begin{eqnarray*}
\Pi[\F_n^c] &\leq& \sum_{h=1}^\infty Q(K=k) \pi_k\left[ \min_{j\leq k} x_j \leq M_n^{-1}\right] \\ 
	    &\leq& \left( \sum_{h=1}^\infty kQ(K=k) \right) \alpha([0,M_n^{-1}]) \\ 
	    &\leq& C' e^{-cn^{1/3}\log(n)^{2/3}}. \\
\end{eqnarray*}

\end{proof}

%
%
%
%
%
%
%
%

\section{Adaptation of Theorem 4 of \citet{RouRivBVM} } 
\label{app:B}
%
%
%
%
%
This Theorem is a slight modification of Theorem 2.9 of \citet{ggv00}. The main deference lies in the handling of the denominator $D_n$ in 
$$
\Pi(f : d(f_0,f) \geq J_{0,n} \epsilon_n | \Xn) = { \int_{d(f,f_0) \geq J_{0,n} \epsilon_n} \prod_{i=1}^n  \frac{f(X_i)}{f_0(X_i)} d\Pi(f) \over\int \prod_{i=1}^n  \frac{f(X_i)}{f_0(X_i)} d\pi(f)   } = \frac{N_n}{D_n},
$$
as in general, it require a lower bound on the prior mass of Kullback Leibler neighborhood of $f_0$. Here we prove that under condition \eqref{eq:cond_kull} we have for some constants $c,C>0$
$$
P_0^n(D_n < c e^{-Cn\epsilon_n^2}) = o(1). 
$$
Let $l_n(f)$ be the log likelihood associated with $f$ and define $\Omega_n = \{ (f,\Xn), l_n(f) - l_n(f_0) > -C_1n\epsilon_n^2\}$ for some constant $C_1>0$. Define also $A_n = \{ \Xn, \forall i X_i \leq \theta_n \}$. We thus have 
$$
D_n \geq e^{-C_1 n \epsilon_n^2} \int_{S_n(\epsilon_n,\theta_n)} \I_{\Omega_n} d\Pi(f) = e^{-C_1n \epsilon_n^2} \Pi(S_n(\epsilon_n, \theta_n) \cap \Omega_n) .
$$
Note that given \eqref{eq:cond_kull} we have that there exists $\rho >0$ such that for $n$ large enough  $e^{-C_2n\epsilon_n^2}\Pi(S_n(\epsilon_n,\theta_n) > \rho$. We now write 
\begin{align*}
	P_0^n (D_n < e^{-Cn\epsilon_n^2}) &\leq P_0^n\left(e^{(C-C_1)n\epsilon_n^2} \Pi(S_n(\epsilon_n,\theta_n)\cap \Omega_n)< c \right) \\
	& \leq P_0^n \left( e^{(C-C_1-C_2)n\epsilon_n^2} \Pi(S_n(\epsilon_n,\theta_n) \cap \Omega_n < \frac{c}{\rho} \Pi(S_n(\epsilon_n,\theta_n) \right) \\ 
	& \leq P_0^n \left( \Pi(S_n(\epsilon_n,\theta_n) \cap \Omega_n^c) > \left( 1 - e^{-(C - C_1 - C_2)n \epsilon_n^2}\frac{c}{\rho}\right)\Pi(S_n(\epsilon_n,\theta_n))  \right) \\
	& \leq { 2 \int_{S_n(\epsilon_n,\theta_n)} P_0^n(\Omega_n^c) d\Pi(f) \over \Pi(S_n(\epsilon_n,\theta_n))}.
\end{align*}
For all $f \in S_n(\epsilon_n,\theta_n)$ we compute 
\begin{align*}
	m_n &= \E_0^n(l_n(f_0) - l_n(f) \I_{A_n}) \\
	&= nF_0(\theta_n)^{n-1} \int_0^{\theta_n} f_{0} \log\left( \frac{f_0(x)}{f(x)} \right)dx \\ 
	&= nF_0(\theta_n)^n\left( KL(f_{0,n},f_n) + \log\left( {F_0(\theta_n) \over F(\theta_n)} \right) \right)\\
	&\leq C_3n\epsilon_n^2,
\end{align*}
and 
\begin{align*}
	P_0^n(\Omega_n^c) &= P_0^n(l_n(f) - l_n(f_0) < -C_1n\epsilon_n^2) \\ 
	&= P_0^n(\{l_n(f) - l_n(f_0) < -C_1 n \epsilon_n^2\} \cap A_n) + o(1) \\ 
	& \leq P_0^n(\{l_n(f_0) - l_n(f) - m_n > (C_1-C_3) n \epsilon_n^2 \} \cap A_n) + o(1) \\
	& \leq { \E_0^n \left(\{l_n(f_0) - l_n(f) - m_n  \} \I_{ A_n}  \right)^2 \over (C_1-C_3)^2(n\epsilon_n^2)^2 } + o(1) .
\end{align*}
We then compute for $C_5$ and $C_6$ some fixed constants
\begin{align*}
	v_n &= \E_0^n\left(\{l_n(f_0) - l_n(f) - m_n  \} \I_{ A_n}  \right)^2  \\
	&= (F_0(\theta_n))^{n-1} \Bigg( n \int_0^{\theta_n} f_{0} \log^2\left(\frac{ f_0(x)}{f(x)} \right) dx
	 + n(n-1) \left( \int_0^{\theta_n} f_{0,n} \log\left(\frac{ f_0(x)}{f(x)} \right)dx\right)^2-m_n^2 \Bigg) \\ 
	&= (F_0(\theta_n))^{n-1} \Bigg( n \int_0^{\theta_n} f_{0} \log^2\left(\frac{ f_0(x)}{f(x)} \right) dx
	 + \frac{n-1}{n} F_0(\theta_n)^{-2n+2}m_n^2-m_n^2 \Bigg) \\ 
	&\leq nF_0(\theta_n)^n  \int_0^{\theta_n}f_{0,n} \log^2\left( \frac{f_0(x)}{f(x)} \right) dx + \frac{n-1}{n} m_n^2 F_0(\theta_n)^{n-1} (F_0(\theta_n)^{-2n+2} -1) \\ 
	& \leq  C_5n\epsilon_n^2 + C_6 (n\epsilon_n^2)^2 \epsilon_n .
\end{align*}
We finally obtain that for all $f \in S_n(\epsilon_n,\theta_n)$, $P_0^n(\Omega_n^c) = o(1)$. We end the proof using similar arguments as in \citet{ggv00}. 
\end{document}